\documentclass[12pt]{amsart}
\usepackage[latin1]{inputenc}
\usepackage{indentfirst}
\usepackage{amsfonts}
\usepackage{graphicx}
\usepackage{amsmath}
\usepackage{amssymb}
\usepackage{latexsym}
\usepackage{amscd}
\usepackage{xcolor}
\usepackage{float}
\usepackage[pagebackref]{hyperref}

\newcommand{\F}{\mathbb {F}}

\newtheorem{theorem}{Theorem}[section]
\newtheorem{definition}[theorem]{Definition}
\newtheorem{lemma}[theorem]{Lemma}
\newtheorem{corollary}[theorem]{Corollary}
\newtheorem{proposition}[theorem]{Proposition}
\newtheorem{remark}[theorem]{Remark}

\def \ord {\rm{ord}\,}

\begin{document}

\title[$r$-primitive and $k$-normal elements]{About $r$- primitive and $k$-normal elements
in finite fields}

\author{Cícero Carvalho}
\author{Josimar J.R. Aguirre} 
\author{Victor G.L. Neumann}
\address{Faculdade de Matemática, Universidade Federal de Uberlandia}
\email{cicero@ufu.br, josimar.mat@ufu.br, victor.neumann@ufu.br}

\maketitle


\vspace{8ex}

\begin{abstract}
In 2013,
Huczynska, Mullen,
Panario and Thomson
introduced the concept of $k$-normal elements: an element $\alpha \in \F_{q^n}$ is $k$-normal over $\F_q$ if the greatest common divisor of the polynomials $g_{\alpha}(x)= \alpha x^{n-1}+\alpha^qx^{n-2}+\ldots +\alpha^{q^{n-2}}x+\alpha^{q^{n-1}}$ and $x^n-1$ in $\F_{q^n}[x]$ has degree $k$,  generalizing the concept of normal elements (normal in the usual sense is $0$-normal).
In this paper we discuss the existence of $r$-primitive,
$k$-normal elements in $\F_{q^n}$ over $\F_{q}$, where
an element $\alpha \in \mathbb{F}_{q^n}^*$ is $r$-primitive if its multiplicative order is $\frac{q^n-1}{r}$. We provide many general results about
the existence of this class of elements and we work a
numerical example over finite fields of characteristic $11$.
\end{abstract}

\section{Introduction}

Let $\mathbb{F}_{q^n}$ be a finite field with $q^n$ elements, where $q$ 
is a prime power and $n$ is a positive integer.
An element $\alpha \in \mathbb{F}_{q^n}^*$ is primitive if $\alpha$ generates
the cyclic multiplicative group $\mathbb{F}_{q^n}^*$. Also, $\alpha \in \mathbb{F}_{q^n}$ is normal over $\mathbb{F}_{q}$ if the set
$B_{\alpha}=\{ \alpha^{q^i} \mid 0 \leq i \leq n-1 \}$ spans $\mathbb{F}_{q^n}$
as an $\mathbb{F}_q$-vector space, in this case  
$B_{\alpha}$ is called a normal basis.
The \textit{Primitive Normal Basis Theorem} states that for any extension field $\mathbb{F}_{q^n}$ of $\mathbb{F}_q$, there exists a basis composed of primitive normal elements; this result was first proved by Lenstra and Schoof \cite{lenstra} using a combination of character sums, sieving results and computer search.

A well-known criterion for checking whether an element generates a normal basis is given by the
following theorem.

\begin{theorem}(\cite[Theorem 2.39]{LN})
For $\alpha \in \mathbb{F}_{q^n}$, $\{\alpha, \alpha^q, \ldots, \alpha^{q^{n-1}} \}$ is a normal basis of $\mathbb{F}_{q^{n}}$ over $\mathbb{F}_q$ if and only if the polynomials $x^n-1$ and $\alpha x^{n-1}+\alpha^qx^{n-2}+\ldots +\alpha^{q^{n-2}}x+\alpha^{q^{n-1}}$ in $\mathbb{F}_{q^n}[x]$ are relatively prime.
\end{theorem}

Inspired by the latest result, the notion of $k$-normal elements is a generalization of the notion of normal elements that was introduced by Huczynska et al. (\cite{knormal}).

\begin{definition}
Let $\alpha \in \mathbb{F}_{q^n}$ and let $g_{\alpha}(x) = \sum_{i=0}^{n-1} \alpha^{q^i}x^{n-1-i} \in \mathbb{F}_{q^n}[x]$. If $\gcd(x^n-1,g_{\alpha}(x))$ over $\mathbb{F}_{q^n}$ has degree $k$ (where $0 \leq k \leq n-1$), then $\alpha$ is a $k$-normal element of $\mathbb{F}_{q^n}$ over $\mathbb{F}_q$.
\end{definition}
Clearly, normal elements are $0$-normal elements. This definition opened a new line of research on the existence of primitive $k$-normal elements (see \cite{lucas}, \cite{lucas1}, \cite{AN} for some results). 


\begin{definition}
An element $\alpha \in \mathbb{F}_{q^n}^*$ is called $r$-primitive if $\alpha$ has order
$(q^n-1)/r$ of $\mathbb{F}_{q^n}^*$, where $r \mid q^n-1$.
\end{definition}

So, primitive elements in the usual sense are $1$-primitive elements. A natural generalization of the problem about the existence of primitive $k$-normal elements, would be to study the existence of $r$-primitive, $k$-normal elements in finite fields. Besides that, it would help to answer another of the questions posed by Huczynska et al. (see \cite[Problem 6.4]{knormal}):  Determine
the existence of high-order $k$-normal elements $\alpha \in \F_{q^n}$ over $\F_q$ , where ``high order"
means $\ord(\alpha)= N$, with $N$ a large positive divisor of $q^n-1$.

In this paper we are going to present some results on the existence of these elements, 
showing that the existing results about the existence of primitive, $k$-normal elements are 
particular cases of the cases treated here. 
In Section $2$, we provide the background material that is used
throughout the paper. In Section $3$, we present two general conditions for the existence
of $r$-primitive $k$-normal elements in $\mathbb{F}_{q^n}$ over $\mathbb{F}_q$, 
as well as some weaker conditions for some particular cases. In the last section  we provide some numerical examples over finite
fields of characteristic $11$.

\section{Preliminaries}
In this section, we present some definitions and results that will be used in
this paper.
We refer the reader to \cite{LN} 
for
basic results on finite fields.

For a positive integer $n$, $\varphi(n)$ denotes the Euler totient function and $\mu(n)$ denotes
the M\"obius function.
\begin{definition}
\begin{enumerate}
\item[(a)] Let $f(x)\in \mathbb{F}_{q}[x]$. The Euler totient function for polynomials over $\mathbb{F}_q$ is given by
$$
\Phi_q(f)= \left| \left( \dfrac{\mathbb{F}_q[x]}{\langle f \rangle} \right)^{*} \right|,
$$
where $\langle f \rangle$ is the ideal generated by $f(x)$ in $\mathbb{F}_q[x]$.
\item[(b)] If $t$ is a positive integer or a monic polynomial over $\mathbb{F}_q$, W(t) denotes the number of squares-free or monic square-free divisors of $t$ respectively.
\item[(c)] If $f(x)\in \mathbb{F}_{q}[x]$ is a monic polynomial, the polynomial M\"obius function $\mu_q$ is given by $\mu_q(f)=0$ if $f$ is not square-free, and $\mu_q(f)=(-1)^r$ if $f$ is a product of $r$ distinct monic irreducible factors over $\mathbb{F}_q$.
\end{enumerate}
\end{definition}

A multiplicative character $\eta$ 
of $\F_{q^n}^*$ is a
group homomorphism of $\F_{q^n}^*$ to $\mathbb{C}^*$.
The group of mutliplicative characters $\widehat{\mathbb{F}}_{q^n}^*$ becomes 
a $\mathbb{Z}$-module by defining
$\eta^r(\alpha)=\eta(\alpha^r)$ for $\eta \in \widehat{\mathbb{F}}_{q^n}^*$,
$\alpha \in \mathbb{F}_{q^n}^*$ and $r \in \mathbb{Z}$.
The order of a multiplicative character $\eta$ is the least positive integer $d$ such that
$\eta (\alpha)^d=1$ for any $\alpha \in \F_{q^n}^*$.

Let $m$ be a divisor of $q^n-1$. We say that an element 
$\alpha \in \mathbb{F}_{q^n}^*$ is $m$-free
if for every $d \mid m$ such that $d\neq 1$ there is no element
$\beta \in \F_{q^n}$ satisfying 
$\alpha= \beta^d$.
Following e.g. \cite[Theorem 13.4.4]{galois} we have that
for any $\alpha \in \F_{q^n}^*$ we get
$$w_m(\alpha) = 
\theta(m) \int_{d|m} \eta_d(\alpha)
=\left\{
\begin{array}{ll}
	1, \quad & \text{if } \alpha \text{ is } m\text{-free}, \\
	0,            & \text{otherwise,}
\end{array}
\right.
$$
where $\theta(m)=\frac{\varphi(m)}{m}$,
$\displaystyle\int_{d|m} \eta_d$ denotes the sum
$\displaystyle \sum_{d|m} \frac{\mu(d)}{\varphi(d)} \sum_{(d)} \eta_d$,
$\eta_d$ is a multiplicative character of $\F_{q^n}$, and the sum
$\displaystyle \sum_{(d)} \eta_d$ runs over all the multiplicative characters of order $d$.

The additive group $\mathbb{F}_{q^n}$ is an $\mathbb{F}_q[x]$-module where the action
is given by $f \circ \alpha =\displaystyle  \sum_{i=0}^r a_i \alpha^{q^i}$,
for any $f=\displaystyle \sum_{i=0}^r a_ix^i\in \mathbb{F}_q[x]$
and $\alpha \in \mathbb{F}_{q^n}$.
An element $\alpha \in \F_{q^n}$ has $\F_q$-order $h\in \mathbb{F}_q[x]$
if $h$ is the monic polynomial of lowest degree such that $h \circ \alpha=0$.
The $\mathbb{F}_q$-order of $\alpha$ will be denoted by $\mathrm{Ord} (\alpha)$, 
and clearly the $\F_q$-order of an element $\alpha \in \F_{q^n}$
divides $x^n-1$. 
An additive character $\chi$ 
of $\F_{q^n}$ is a
group homomorphism of $\mathbb{F}_{q^n}$ to $\mathbb{C}^*$.
The
group of additive characters $\widehat{\mathbb{F}}_{q^n}$ becomes 
an $\mathbb{F}_q[x]$-module by defining
$f\circ \chi (\alpha)=\chi(f \circ \alpha)$ for $\chi \in \widehat{\mathbb{F}}_{q^n}$,
$\alpha \in \mathbb{F}_{q^n}$ and $f \in \mathbb{F}_q[x]$.
An additive character $\chi$ has $\F_q$-order $h \in \mathbb{F}_q[x]$
if $h$ is the monic polynomial of smallest degree such
that $h \circ \chi$
is the trivial additive character.
The $\mathbb{F}_q$-order of $\chi$ will be denoted by $\mathrm{Ord} (\chi)$.

Let $g\in \mathbb{F}_q[x]$ be a  divisor of $x^n-1$. We say that an element
$\alpha \in \mathbb{F}_{q^n}$ is $g$-free
if for every polynomial $h \in \mathbb{F}_q[x]$  such that $h \mid g$ and $h\neq 1$, there is no element
$\beta \in \F_{q^n}$ satisfying $\alpha= h \circ \beta$.
As in the multiplicative case, from e.g. \cite[Theorem 13.4.4]{galois} we have that
for any $\alpha \in \F_{q^n}$ we get
$$
\Omega_g(\alpha)=
\Theta(g) \int_{h|g} \chi_h(\alpha)
=\left\{
\begin{array}{ll}
	1, \quad & \text{if } \alpha \text{ is } g\text{-free}, \\
	0,            & \text{otherwise,}
\end{array}
\right.
$$
where 
$\Theta(g)= \frac{\Phi_q(g)}{q^{\deg(g)}}$,
$\displaystyle\int_{h|g} \chi_h$ denotes the sum
$\displaystyle \sum_{h|g} \frac{\mu_q(h)}{\Phi_q(h)} \sum_{(h)} \chi_h$,
$\displaystyle \sum_{h|g}$ runs over all the monic divisors $h\in \mathbb{F}_q[x]$ of $g$,
$\chi_h$ is an additive character of $\F_{q^n}$, and the sum
$\displaystyle \sum_{(h)} \chi_h$ runs over all additive characters of $\mathbb{F}_q$-order $h$.
It is known that there exist $\Phi_q(h)$ of those characters.



One may check that an element $\alpha \in \mathbb{F}_{q^n}^*$ is primitive if and only if $\alpha$ is $(q^n-1)$-free and $\alpha \in \mathbb{F}_{q^n}$ is normal if and only if $\alpha$ is $(x^n-1)$-free. 

%
\begin{remark}\label{construct-k-r}
In \cite{lucas} Reis gives a method to construct $k$-normal elements: let $\beta \in \mathbb{F}_{q^n}$ be a normal element and $f\in \mathbb{F}_q[x]$ be a divisor of $x^n-1$ of degree $k$, then $\alpha = f \circ \beta$ is $k$-normal
(see \cite[Lemma 3.1]{lucas}). In the same way, if $\beta \in \mathbb{F}_{q^n}$ is a primitive element, then $\beta^r$ is $r$-primitive
for any divisor $r$ of $q^n-1$.
\end{remark}

We have that $\mathbb{F}_{q^n}^*$ and $\widehat{\mathbb{F}}_{q^n}^*$
are isomorphic as $\mathbb{Z}$-modules, and
$\mathbb{F}_{q^n}$ and $\widehat{\mathbb{F}}_{q^n}$
are isomorphic as $\mathbb{F}_q[x]$-modules (see \cite[Theorem 13.4.1]{galois}).


We will also need the following definition to characterize when a k-normal element can also be r-primitive.
\begin{definition}\label{char0}
For any $\alpha \in \mathbb{F}_{q^n}$ we define the following character sum:
$$
I_0(\alpha) = \frac{1}{q^n} \sum_{\psi \in \widehat{\mathbb{F}}_{q^n}} \psi(\alpha).
$$
Note that $I_0(\alpha)=1$ if $\alpha = 0$, and $I_0(\alpha)=0$ otherwise by the character orthogonality property.
\end{definition}

To finish this section, we present some estimates that are used in the next sections.


\begin{lemma}\label{case-a}
(\cite[Theorem 5.6]{Fu}) Let $r \in \mathbb{N}$ be a divisor of $q^n-1$,
let $\eta$ be a multiplicative character and let $\psi$ be a non-trivial  additive character.
Then
$$
\left| \sum_{\alpha \in \mathbb{F}_{q^n}^*} \eta(\alpha) \psi(\alpha^r)\right|
\leq r q^{n/2}.
$$
\end{lemma}

\begin{lemma}\label{case-b}
Let $f \in \mathbb{F}_q[x]$ be a divisor of $x^n-1$ of degree $k$ and
let $\chi$ and $\psi$ be additive characters. Then
$$
\sum_{\beta \in \mathbb{F}_{q^n}}
\chi(\beta)\psi(f \circ \beta)^{-1} =
\left\{
\begin{array}{ll}
	q^n & \text{if } \chi = f \circ \psi, \\
	0 & \text{if } \chi \neq f \circ \psi .
\end{array}
\right.
$$
Furthermore, for a given additive character $\chi$, 
the set
$\{\psi \in  \widehat{\mathbb{F}}_{q^n} \mid \chi = f \circ \psi\}$
has $q^k$ elements if 
$\mathrm{Ord}(\chi) \mid \frac{x^n-1}{f}$, and it is an empty set
if $\mathrm{Ord}(\chi) \nmid \frac{x^n-1}{f}$.
\end{lemma}
\begin{proof}
From \cite[Theorem 5.4]{LN}, the sum is zero if and only if $\chi \neq f \circ \psi$,
since for any $\beta \in \mathbb{F}_{q^n}$ we have
$\chi(\beta)\psi(f \circ \beta)^{-1}=(\chi \cdot (f \circ \psi)^{-1})(\beta)$, and
if $\chi = f \circ \psi$ the sum is $q^n$.

For any  additive character $\chi$, we have that $\mathrm{Ord}(\chi) \mid \frac{x^n-1}{f}$ if and only if
$\frac{x^n-1}{f} \circ \chi$ is the trivial character.
From \cite[Theorem 13.4.1]{galois}, $f$ and $\frac{x^n-1}{f}$ define linear endomorphisms of
$\widehat{\mathbb{F}}_{q^n}$ like they do in  $\mathbb{F}_{q^n}$, so from duality and 
\cite[Lemma 2.5]{AN}
we have that
$\frac{x^n-1}{f} \circ \chi$ is the trivial character if and only
if there exists an additive character $\psi$ such that $\chi = f\circ \psi$. 
This proves that $\{\psi \in  \widehat{\mathbb{F}}_{q^n} \mid \chi = f \circ \psi\}\neq \emptyset$
if and only if $\mathrm{Ord}(\chi) \mid \frac{x^n-1}{f}$.

Denote by $\hat{f}$ the linear endomorphism of $\widehat{\mathbb{F}}_{q^n}$ defined by $f$. From duality
$\ker \hat{f}$ has $q^k$ elements, since $\ker f$ has $q^k$ elements (see the proof of \cite[Lemma 2.5]{AN}). So, if 
$\chi \in \mathop{\mathrm{im}} \hat{f}$, the set
$\hat{f}^{-1}(\chi) = \{\psi \in  \widehat{\mathbb{F}}_{q^n} \mid \chi = f \circ \psi\}$
has $q^k$ elements, since
$\hat{f}^{-1}(\chi)$ is a coset of 
$\ker \hat{f}$.
\end{proof}

%
%

\section{General results}\label{sectiongen}
We are interested in finding conditions for the existence of $r$-primitive $k$-normal elements. For this, the following definition plays an important role.

\begin{definition}\label{def-NrfmT}
Let $f,g \in \F_q[x]$ be monic divisors of $x^n-1$, with $\deg f =k$,
and let $m,r \in \mathbb{N}$ be divisors
of $q^n-1$. 
We define
$$
N_{r,f}(m,g) = \sum_{\alpha \in \mathbb{F}_{q^n}^*}
\sum_{\beta \in \mathbb{F}_{q^n}}
w_m(\alpha) \Omega_g(\beta) I_0(\alpha^r - f \circ \beta).
$$
From the definition of $w_m$, $\Omega_g$ and Definition \ref{char0},
$N_{r,f}(m,g)$
counts the number of pairs $(\alpha,\beta) \in \mathbb{F}_{q^n}^* \times \mathbb{F}_{q^n}$
such that $\alpha$ is $m$-free, $\beta$ is $g$-free and $\alpha^r=f \circ \beta$.
In particular, if $N_{r,f}(q^n-1,x^n-1)>0$ then there exists a pair 
$(\alpha,\beta) \in \mathbb{F}_{q^n}^* \times \mathbb{F}_{q^n}$
such that 
$\alpha$ is primitive, $\beta$ is normal and
$\alpha^r=f \circ \beta$. From Remark \ref{construct-k-r},
$\alpha^r=f \circ \beta$
is an $r$-primitive and $k$-normal element of $\mathbb{F}_{q^n}$
over $\mathbb{F}_q$.
\end{definition}

Using the last definition, we need to find lower estimates for the sum above, in order to guarantee the positivity of $N_{r,f}(q^n-1,x^n-1)$. We have the following result:

\begin{theorem}\label{principal}
Let $q$ be a prime power and let $n$ be a positive integer.
Let $r \in \mathbb{N}$ be a  divisor
of $q^n-1$  and
let $f \in \F_q[x]$ be a monic divisor of $x^n-1$ of degree $k$.
If 
$q^{\frac{n}{2}-k} \geq r W(q^n-1)W(\frac{x^n-1}{f})$
then 
there exists an $r$-primitive and $k$-normal element of $\mathbb{F}_{q^n}$ over $\mathbb{F}_q$.
\end{theorem}
\begin{proof}
Let $m \in \mathbb{N}$ be a divisor of $q^n-1$
and let $g \in \mathbb{F}_q[x]$ be a divisor of $x^n-1$. We will find
a bound for
$N_{r,f}(m,g)$.
From the definitions of $w_m$, $\Omega_g$ and Definition \ref{char0}, we have
$$
N_{r,f}(m,g)  = \frac{\theta(m)\Theta(g)} {q^n}  
\int_{d \mid m} \int_{h \mid g} \sum_{{\psi \in \widehat{\mathbb{F}}_{q^n}} }
S(\eta_d, \chi_h, \psi),
$$
where
\begin{align*}
S(\eta, \chi, \psi) & = 
\sum_{\alpha \in \mathbb{F}_{q^n}^*}
\sum_{\beta \in \mathbb{F}_{q^n}}
\eta (\alpha) \chi (\beta) \psi(\alpha^r - f \circ \beta) \\
& = 
\sum_{\alpha \in \mathbb{F}_{q^n}^*} \eta(\alpha) \psi(\alpha^r )
\sum_{\beta \in \mathbb{F}_{q^n}}
\chi (\beta) \psi(-f \circ \beta).
\end{align*}
Let $\eta_1$ be the trivial multiplicative character and
let $\chi_1$ be the trivial additive character. Write
$$
N_{r,f}(m,g) = \frac{\theta(m)\Theta(g)} { q^n} 
\left(
S_1 + S_2 + S_3 + S_4
\right),
$$
where
$
S_1=S(\eta_1,\chi_1,\chi_1),$ 
$$
S_2= 
\mathop{\int}_{d \mid m}
\sum_{\substack{\psi \in \widehat{\mathbb{F}}_{q^n} \\ \psi \neq \chi_1}}
S (\eta_d,\chi_1,\psi),
\hspace*{1cm}
S_3=
\mathop{\int}_{\substack{d \mid m \\ h \mid g \\ d\neq 1 \text{ or } h \neq 1}}
S (\eta_d,\chi_h,\chi_1),
$$
and
$$
S_4= 
\mathop{\int}_{d \mid m}
\mathop{\int}_{\substack{h \mid g \\ h \neq 1}}
\sum_{\substack{\psi \in \widehat{\mathbb{F}}_{q^n} \\ \psi \neq \chi_1}}
S (\eta_d,\chi_h,\psi).
$$
We have $S_1=(q^n-1)q^n$, since $\eta_1(\alpha) \chi_1(\alpha^r )=1$ and
$\chi_1 (\beta) \chi_1(-f \circ \beta)=1$ for every $\alpha \in \mathbb{F}_{q^n}^*$
and $\beta \in \mathbb{F}_{q^n}$.

From Lemma \ref{case-b},
for any multiplicative character $\eta_d$, we have
$$
\sum_{\substack{\psi \in \widehat{\mathbb{F}}_{q^n} \\ \psi \neq \chi_1}}
S (\eta_d,\chi_1,\psi) 
  =
\sum_{\substack{\psi \in \ker \hat{f} \\ \psi \neq \chi_1}}
S (\eta_d,\chi_1,\psi)  
= 
(q^k-1) q^n  \sum_{\alpha \in \mathbb{F}_{q^n}^*} \eta_d(\alpha) \psi(\alpha^r ),
$$
where $\ker \hat{f}=\{\psi \in  \widehat{\mathbb{F}}_{q^n} \mid \chi_1 = f \circ \psi\}$.
Now, from Lemma \ref{case-a} and using that
there are $\varphi(d)$ multiplicative characters of order $d$,
we get
$$
|S_2| \leq \sum_{d\mid m} \frac{|\mu(d)|}{\varphi(d)} \sum_{(d)} r (q^k-1) q^{\frac{3n}{2}}
=r (q^k-1) q^{\frac{3n}{2}} W(m).
$$

From \cite[Theorem 5.4]{LN}, we have $S(\eta,\chi,\chi_1)=0$ if $\eta \neq \eta_1$ or $\chi \neq \chi_1$, so
$S_3=0$.

To get a bound of $S_4$ we define $\widetilde{g}=\gcd (g, \frac{x^n-1}{f})$ and, 
for any additive character $\chi$, consider
$\hat{f}^{-1}(\chi)=\{\psi \in  \widehat{\mathbb{F}}_{q^n} \mid \chi = f \circ \psi\}$.
From Lemma \ref{case-b}, we have that $\hat{f}^{-1}(\chi)=\emptyset$ if
$\mathrm{Ord}(\chi) \nmid  \frac{x^n-1}{f}$, so
$$
S_4 = 
\mathop{\int}_{d \mid m}
\mathop{\int}_{\substack{h \mid \widetilde{g} \\ h \neq 1}}
\sum_{\psi \in \hat{f}^{-1}(\chi_h)}
S (\eta_d,\chi_h,\psi).
$$
Again from Lemma \ref{case-b}, for 
any multiplicative character $\eta_d$ of order $d$ and
any additive character $\chi_h$ whose $\mathbb{F}_q$-order $h$ divides $\widetilde{g}$, we have
$$
\sum_{\psi \in \hat{f}^{-1}(\chi_h)}
S (\eta_d,\chi_h,\psi)  
= 
q^{n+k} \sum_{\alpha \in \mathbb{F}_{q^n}^*} \eta_d(\alpha) \psi(\alpha^r ).
$$
So, from Lemma \ref{case-a},
$$
|S_4| \leq \sum_{d\mid m} \sum_{\substack{h \mid \widetilde{g} \\ h \neq 1}}
\frac{|\mu(d) \mu_q(h)|}{\varphi(d) \Phi_q(h)} \sum_{(d)} \sum_{(h)} r q^{\frac{3n}{2}+ k}
=r q^{\frac{3n}{2}+k} W(m) (W(\widetilde{g})-1).
$$
Therefore, we conclude that
\begin{eqnarray*} 
N_{r,f}(m,g) & \geq &\dfrac{\theta(m) \Theta(g)}{q^n}
		\left(
		(q^n-1)q^n- r (q^k-1) q^{\frac{3n}{2}} W(m)  \right.  \\
		& & \left. -  r q^{\frac{3n}{2}+k} W(m) (W(\widetilde{g})-1)
		\right) \\
		& > & \theta(m) \Theta(g) \left( q^n - r q^{\frac{n}{2}+k} W(m) W(\widetilde{g}) \right),
\end{eqnarray*}
since $rq^{\frac{n}{2}}W(m) - 1>0$. Thus, if
$q^{\frac{n}{2}-k} \geq rW(m) W(\widetilde{g})$ then
$N_{r,f}(m,g)>0$. In particular,
If 
$q^{\frac{n}{2}-k} \geq r W(q^n-1)W(\frac{x^n-1}{f})$
then 
there exists an $r$-primitive and $k$-normal element of $\mathbb{F}_{q^n}$ over $\mathbb{F}_q$.
\end{proof}

\begin{remark}
Notice that the last result, with $r=k=1$ and $f=x-1$, generalizes
previous results on the existence of primitive $1$-normal elements (see \cite[Corollary 5.8]{knormal}).
Also, this is a stronger condition for the existence of primitive, $k$-normal elements 
than the result given in \cite[Theorem 3.3]{lucas}.
\end{remark}

The sieving technique from the next two results is similar to others which have appeared in  
previous works about primitive or normal elements.

\begin{lemma}\label{lemmasieve}
Let $q$ be a prime power and let $n$ be a positive integer.
Let $r \in \mathbb{N}$ be a  divisor of $q^n-1$  and
let $f \in \F_q[x]$ be a divisor of $x^n-1$ of degree $k$.
Let $\ell$ be a divisor of $q^n-1$ and let $\{p_1,...,p_v\}$ be the set of 
all 
primes which divide $q^n-1$,
but do not divide $\ell$.
Also let $g \in \F_q[x]$ be a divisor of $x^n -1$ and $\{P_1,...,P_s\} \subset \F_q[x]$
be the set of all monic irreducible polynomials which
divide $x^n -1$, but do not divide $g$. Then
	\begin{equation}\label{sieve}
		N_{r,f}(q^n-1,x^n-1)\geq 
		\sum_{i=1}^{v}N_{r,f}(p_i \ell,g)
		+
		\sum_{i=1}^{s}N_{r,f}(\ell, P_i g)
		- (v+s-1)N_{r,f}(\ell,g). 
	\end{equation}
\end{lemma}
\begin{proof}
The left side of \eqref{sieve} counts every pair $(\alpha,\beta) \in \mathbb{F}_{q^n}^* \times \mathbb{F}_{q^n}$
for which $\alpha$ is primitive, $\beta$ is normal and $\alpha^r=f \circ \beta$.
Observe that for such a pair $(\alpha,\beta)$ we also have that
$\alpha$ is  $p_i \ell$-free and $\ell$-free, and $\beta$ is $P_i g$-free and $g$-free, so $(\alpha,\beta)$ is  counted $v+s - (v+s-1)=1$
times on the right side of \eqref{sieve}.
For any other pair $(\alpha,\beta) \in \mathbb{F}_{q^n}^* \times \mathbb{F}_{q^n}$, we have that either
$\alpha^r \neq f \circ \beta$, or
$\alpha$
is not $p_i \ell$-free for some $i \in \{ 1,\ldots , v\}$
or $\beta$ is not $P_i g$-free for some
$i \in \{ 1,\ldots , s\}$, so the right side of \eqref{sieve} is at most zero.
\end{proof}

\begin{proposition}\label{sieve-prop}
Let $q$ be a prime power and let $n$ be a positive integer.
Let $r \in \mathbb{N}$ be a  divisor of $q^n-1$  and
let $f \in \F_q[x]$ be a divisor of $x^n-1$ of degree $k$.
Let $\ell$ be a divisor of $q^n-1$ and let $\{p_1,...,p_v\}$ be the set of 
all 
primes which divide $q^n-1$,
but do not divide $\ell$.
Also, let $g \in \F_q[x]$ be a divisor of $x^n -1$ and $\{P_1,...,P_s\} \subset \F_q[x]$
be the set of all monic irreducible polynomials which
divide $x^n -1$ but do not divide $g$,
and if necessary reorder the set $\{ P_1, \ldots, P_s \}$ so that $\{ P_1, \ldots, P_{s'}\}$ is the subset of polynomials which divide $x^n - 1/f$ .
Suppose that 
$\delta=1-\sum_{i=1}^{v}\frac{1}{p_i}-\sum_{i=1}^{s}\frac{1}{q^{\deg P_i}}>0$
and let  
$\delta'=1 - \sum_{i=s'+1}^{s} \frac{1}{q^{\deg P_i}}$ and 
$\Delta=2+ \frac{v+s'-\delta'}{\delta}$. 
If
$q^{\frac{n}{2}-k} \geq r W(\ell) W(\gcd(g,\frac{x^n-1}{f}))  \Delta$
then $N_{r,f}(q^n-1,x^n-1)>0$.
\end{proposition}
\begin{proof}
We can rewrite inequality \eqref{sieve} in the form
\begin{align*}
N_{r,f}(q^n-1,x^n-1) & \geq \sum_{i=1}^v \Big[N_{r,f}(p_i \ell ,g) - \theta(p_i) N_{v,f}(\ell,g) \Big]\\
		& + \sum_{i=1}^s \Big[ N_{r,f}(\ell, P_i g) - \Theta(P_i)N_f(\ell,g) \Big]  + \delta N_{r,f}(\ell,g).
\end{align*}
Let $i \in \{1,\ldots,s\}$. From the definitions of $w_m$, $\Omega_g$, Definition \ref{char0} and 
Definition \ref{def-NrfmT},
taking into account that
$\theta$ is a multiplicative function,
we get
\begin{eqnarray*}
N_{r,f}(p_i \ell,g)  & = & \frac{\theta(p_i) \theta(\ell)\Theta(g)} {q^n}  
\mathop{\int}_{d \mid p_i \ell} \mathop{\int}_{h \mid g} \sum_{{\psi \in \widehat{\mathbb{F}}_{q^n}} }
S(\eta_d, \chi_h, \psi) \\
& = & \theta(p_i) N_{r,f}(\ell,g) +
\frac{\theta(p_i) \theta(\ell)\Theta(g)} {q^n}  
\mathop{\int}_{\substack{d \mid p_i \ell \\ p_i \mid d}} \mathop{\int}_{h \mid g} \sum_{{\psi \in \widehat{\mathbb{F}}_{q^n}} }
S(\eta_d, \chi_h, \psi) .
\end{eqnarray*}
Now, from Lemma \ref{case-b},
denoting
$\widetilde{g}=\gcd(g,\frac{x^n-1}{f})$, we have
$$
\mathop{\int}_{\substack{d \mid p_i \ell \\ p_i \mid d}} \mathop{\int}_{h \mid g} \sum_{{\psi \in \widehat{\mathbb{F}}_{q^n}} }
S(\eta_d, \chi_h, \psi)=
q^n \mathop{\int}_{\substack{d \mid p_i \ell \\ p_i \mid d}} 
\mathop{\int}_{h \mid \widetilde{g}}
\sum_{{\psi \in \hat{f}^{-1}(\chi_h)}}
\sum_{\alpha \in \mathbb{F}_{q^n}^*} \eta(\alpha) \psi(\alpha^r ).
$$
From Lemma \ref{case-a}, using that the second summation on the right side of the above equality is zero if $\psi=\chi_1$, we get
$$
\left|
\mathop{\int}_{\substack{d \mid p_i \ell \\ p_i \mid d}} \mathop{\int}_{h \mid g} \sum_{{\psi \in \widehat{\mathbb{F}}_{q^n}} }
S(\eta_d, \chi_h, \psi)
\right| < 
 r q^{\frac{3n}{2}+k} W(\ell) W(\widetilde{g}),
$$
so
$\left|
N_{r,f}(p_i \ell ,g) - \theta(p_i) N_{r,f}(\ell,g)
\right|<
\theta(p_i) \theta(\ell)\Theta(g) r q^{\frac{n}{2}+k} W(\ell) W(\widetilde{g})$.

Once again, from \cite[Theorem 13.4.4.]{galois}, Definition \ref{char0} and 
Definition \ref{def-NrfmT}, taking into account that
$\Theta$ is a multiplicative function,
we have
\begin{eqnarray*}
	N_{r,f}(\ell,P_i g)  & = & \frac{\Theta(P_i) \theta(\ell)\Theta(g)} {q^n}  
	\mathop{\int}_{d \mid \ell} \mathop{\int}_{h \mid P_i g} \sum_{{\psi \in \widehat{\mathbb{F}}_{q^n}} }
	S(\eta_d, \chi_h, \psi) \\
	& = & \Theta(P_i) N_{r,f}(\ell,g) +
	\frac{\Theta(P_i) \theta(\ell)\Theta(g)} {q^n}  
	\mathop{\int}_{d \mid \ell} \mathop{\int}_{\substack{h \mid P_i g \\ P_i \mid h}}
	 \sum_{{\psi \in \widehat{\mathbb{F}}_{q^n}} }
	S(\eta_d, \chi_h, \psi) .
\end{eqnarray*}
Now, from Lemma \ref{case-b}, 
denoting
$\widetilde{g}=\gcd(g,\frac{x^n-1}{f})$
and
$\widetilde{P_i}=\gcd(P_i,\frac{x^n-1}{f})$, we have
$$
\mathop{\int}_{d \mid  \ell} 
\mathop{\int}_{\substack{h \mid P_i g \\ P_i \mid h}}
\sum_{{\psi \in \widehat{\mathbb{F}}_{q^n}} }
S(\eta_d, \chi_h, \psi)=
q^n \mathop{\int}_{d \mid  \ell} 
\mathop{\int}_{\substack{h \mid \widetilde{P_i} \widetilde{g} \\ P_i \mid h}}
\sum_{{\psi \in \hat{f}^{-1}(\chi_h)}}
\sum_{\alpha \in \mathbb{F}_{q^n}^*} \eta(\alpha) \psi(\alpha^r ).
$$
Note that if $i> s'$  then the right side of the above equality is zero, since $P_i \nmid \frac{x^n-1}{f}$. Suppose that
$1 \leq i \leq s'$.
From Lemma \ref{case-a}, we get
$$
\left|
\mathop{\int}_{d \mid  \ell} 
\mathop{\int}_{\substack{h \mid P_i g \\ P_i \mid h}}
\sum_{{\psi \in \widehat{\mathbb{F}}_{q^n}} }
S(\eta_d, \chi_h, \psi)
\right| \leq
r q^{\frac{3n}{2}+k} W(\ell) W(\widetilde{g}),
$$
so
$\left|
N_{r,f}(\ell ,P_i g) - \Theta(P_i) N_{r,f}(\ell,g)
\right|\leq 
\Theta(P_i) \theta(\ell)\Theta(g) r q^{\frac{n}{2}+k} W(\ell) W(\widetilde{g})$.


Combining all the inequalities above we obtain
	\begin{align*}
		N_{r,f}(q^n-1,x^n-1) & \geq \delta N_{r,f}(\ell,g) - \\
		&  \theta(\ell) \Theta(g) W(\ell)W(\widetilde{g}) r q^{n/2+k} \left( \sum_{i=1}^v \theta(p_i) + \sum_{i=1}^{s'}\Theta(P_i) \right).
	\end{align*}
Therefore, 
from the proof of Theorem \ref{principal}, we have
\begin{align*}
N_{r,f}(q^n-1,x^n-1) & > \delta \theta(\ell) \Theta(g) \left( q^n- r q^{\frac{n}{2}+k}W(\ell)W(\widetilde{g}) \right) \\ &  \ \  
\ - \theta(\ell) \Theta(g) W(\ell)W(\widetilde{g}) r q^{\frac{n}{2}+k} \left( \sum_{i=1}^v \theta(p_i) + \sum_{i=1}^{s'}\Theta(P_i) \right) \\
& = \delta \theta(\ell)\Theta(g) \Big( q^n - r q^{\frac{n}{2}+k}W(\ell)W(\widetilde{g}) \Delta \Big),
\end{align*}
and we obtain the desired result.
\end{proof}

\begin{proposition}\label{caseall} 
Let $n$ be a positive integer and let $q$ be a prime power.
Let $r,k \in \mathbb{N}$ such that $r$ is a divisor
of $q^n-1$, $k<n/2$,
there exists a degree $k$ factor of $x^n-1$ in $\mathbb{F}_q[x]$ and
$(n-k)^2 \leq q$. 
If $q^{\frac{n}{2}-k} \geq r(n-k+2) W(q^n-1)$, then 
there exists an $r$-primitive $k$-normal element in $\F_{q^n}$.
\end{proposition}
\begin{proof}
Let $f \in \F_q[x]$ be a factor of $x^n-1$ of degree $k$.
We may use Proposition \ref{sieve-prop} with $\ell=q^n-1$ and $g$ a divisor of $x^n-1$ such that
$\gcd(g,\frac{x^n-1}{f})=1$ and any irreducible factor of $x^n-1$ divides $g$ or
$\frac{x^n-1}{f}$.

Let $P_1,\ldots , P_s$ be all the  irreducible polynomials such that
$\mathrm{rad}(\frac{x^n-1}{f})= P_1 \cdot P_2 \cdots P_s$.
Then
$\delta=1-\sum_{i=1}^{s}\frac{1}{q^{\deg P_i}} \geq 1 - \frac{n-k}{q}
\geq 1 - \frac{1}{n-k} = \frac{n-k-1}{n-k} >0$, since
$q \geq (n-k)^2$ and $s \leq n-k$.  We also have that
$$
\Delta=2+\frac{s-1}{\delta}\leq \frac{n-k-1}{\frac{n-k-1}{n-k}}+2 = n-k+2.
$$
This means that $W(\ell) W(\widetilde{g})  \Delta \leq (n-k+2)W(q^n-1)$ and from 
Proposition \ref{sieve-prop} we get
	the desired result.
\end{proof}

The next result will be used in specific cases.
\begin{lemma}\label{casetu}
Let $n$ be a positive integer and let $q$ be a prime power.
Let $r,k \in \mathbb{N}$ such that $r$ is a divisor
of $q^n-1$, $k<n/2$,
there exists a degree $k$ factor of $x^n-1$ in $\mathbb{F}_q[x]$ and
$(n-k)^2 \leq q$. Also, let $t,u$ be positive real numbers such that $t+u > \dfrac{2n}{n-2k}$
and $\delta_{t,u}=1 - S_{t,u} - \frac{1}{n-k}>0$, where $S_{t,u}$ is the 
sum of the inverse of all prime numbers between $2^t$ and $2^{t+u}$.
If
$$
q \geq   \left( r 
\Delta_{t,u} A_{t,u}
\right)^{ \frac{2(t+u)}{(t+u)(n-2k)-2n}},
$$
then 
there exists an $r$-primitive $k$-normal element in $\F_{q^n}$, where
$$
A_{t,u}=\prod_{\substack{\wp < 2^t \\ \wp \text{ is prime}}}
\frac{2}{\sqrt[t+u]{\wp}} , \quad
\Delta_{t,u}=2+\frac{v(t,u)+n-k-1}{\delta_{t,u}}
$$
and $v(t,u)$ is  the number of
all prime numbers between $2^t$ and $2^{t+u}$.
\end{lemma}
\begin{proof}
Let $t,u$ be positive real numbers such that $t+u > \dfrac{2n}{n-2k}$ and
let 
$$
q^n-1=p_1^{a_1} \cdots p_w^{a_{w}} \cdot
	q_1^{b_1} \cdots q_v^{b_{v}}
	$$
be the prime factorization of $q^n-1$ such that $2 \leq p_i \leq 2^t$  or 
	$2^{t+u} \leq p_i$
	for $1 \leq i \leq w$ and
	$2^t < q_i < 2^{t+u}$
	for $1 \leq i \leq v$.
	We use Proposition \ref{sieve-prop}, where
$\ell = p_1^{a_1} \cdots p_w^{a_{w}}$
and 
$g$ is a divisor of $x^n-1$ such that
$\gcd(g,\frac{x^n-1}{f})=1$, and any irreducible factor of $x^n-1$ divides $g$ or
$\frac{x^n-1}{f}$.

Let $P_1,\ldots , P_s$ be all the  irreducible polynomials such that
$\mathrm{rad}(\frac{x^n-1}{f})= P_1 \cdot P_2 \cdots P_s$.
Then
$\delta=1 - \sum_{i=1}^v \frac{1}{q_i} -\sum_{i=1}^{s}\frac{1}{q^{\deg P_i}} \geq 1 - \sum_{i=1}^v \frac{1}{q_i} - \frac{n-k}{q}
\geq 1 - S_{t,u} - \frac{1}{n-k}$,
where $\mathrm{rad} (x^n-1) = Q_1 \cdots Q_s$. If
$\delta_{t,u}>0$ then
$\Delta =2 + \frac{v+s-1}{\delta} \leq 2+\frac{v(t,u)+n-k-1}{\delta_{t,u}}=\Delta_{t,u}$.
Let  $P_t$ be the set of all prime numbers less than $2^t$.
From \cite[Lemma 3.7]{cgnt2} we get $W(\ell) \leq A_{t,u} \ell^{\frac{1}{t+u}} \leq A_{t,u} q^{\frac{n}{t+u}}$.
From Proposition \ref{sieve-prop}, we conclude that a sufficient condition for the existence of 
a $r$-primitive, $k$-normal element in $\F_{q^n}$ is
$q^{\frac{n}{2}-k} \geq r  \Delta_{t,u}  A_{t,u}  q^{\frac{n}{t+u}}$
or, equivalently,
$$
q \geq   \left( r 
	\Delta_{t,u}  A_{t,u}
	\right)^{ \frac{2(t+u)}{(t+u)(n-2k)-2n}}.
$$
\end{proof}

The next result shows that if 
$n$, $k$ and $r$  are positive integers such that $k < n/2$
then there exists a constant $C(n, k, r)$ such that if  
$q \geq C(n,k,r)$,
$r \mid q^n-1$ and
there exists a degree $k$ factor of $x^n-1$ in $\mathbb{F}_q[x]$,
then there exists an $r$-primitive  element in $\F_{q^n}$ which 
is $k$-normal over 
$\F_q$.

\begin{proposition}\label{cota}
Let $n$, $k$ and $r$  be positive integers such that $r \mid q^n-1$, $k < n/2$ and
there exists a degree $k$ factor of $x^n-1$ in $\mathbb{F}_q[x]$.
Let $t$ be a 
real 
number such that $t >  
2n/(n - 2k)$. Then if
\begin{equation}\label{condition}
q \geq \min \{  U_t(n,k,r) , \max  \{(n-k)^2,V_t(n,k,r) \}   \}
\end{equation}
then 
there exists an $r$-primitive element in $\F_{q^n}$ which is $k$-normal over 
$\F_q$, where
\begin{eqnarray*}
U_t(n,k,r) & = & \left( r 2^{n - k} A_t  \right)^{\frac{2t}{t(n - 2k)- 2n}} ,\\
V_t(n,k,r) & = & \left( r (n - k+2) A_t  \right)^{\frac{2t}{t(n - 2k)- 2n}} \quad \text{and}\\
A_t & = & \prod_{\substack{\wp < 2^t \\ \wp \text{ is prime}}}
\frac{2}{\sqrt[t]{\wp}} .
\end{eqnarray*}
\end{proposition}
\begin{proof}
From \cite[Lemma 3.7]{cgnt2} we  get 
\[
W(q^n - 1) \leq A_t q^{n/t}.
\]	
Clearly $W(\frac{x^n-1}{f}) \leq 2^{n - k}$, so from Theorem  \ref{principal} we 
get that if 
\[
q^{\frac{n}{2}-k} \geq  r q^{n/t} A_t 2^{n - k}, 
\]
or equivalently, if
\[
q \geq \left( r 2^{n - k} A_t  \right)^{\frac{2t}{t(n - 2k)- 2n}}
\]
 then 
there exists an $r$-primitive element in $\F_{q^n}$ which is $k$-normal over 
$\F_q$.
Now, if $q\geq (n-k)^2$ then, from Proposition \ref{caseall} and \cite[Lemma 3.7]{cgnt2}, we get that
if
\[
q \geq \left( r (n - k+2) A_t  \right)^{\frac{2t}{t(n - 2k)- 2n}}
\]
then 
there exists an $r$-primitive element in $\F_{q^n}$ which is $k$-normal over 
$\F_q$.
\end{proof}

%
%
%

\section{Numerical example}
We apply Proposition \ref{cota} to study  $3$-primitives elements in $\mathbb{F}_{q^n}$
which are $3$-normal over $\mathbb{F}_q$, where $\mathbb{F}_q$ is a field of characteristic $11$.
\begin{proposition}\label{caser3k3}
Let $n\geq 7$ be a positive integer and let $q$ be a prime power.
For any pair $(q,n)$ in Table \ref{table1}, if $3 \mid q^n-1$ and
$x^n-1$ has a degree $3$ factor in $\mathbb{F}_q[x]$ then
there exists a $3$-primitive element in  $\F_{q^n}$ which is $3$-normal over $\F_{q}$.
\begin{table}[h]
	\centering
	\begin{tabular}{c|c}
		$t$ or $(t,u)$ & $(q,n)$ \\
		\hline 
		$t=7.5$ & $q \geq 11$ and $n \geq 70$ \\
		\hline 
		$t=7$ & $q \geq 16$ and $n \geq 44$ \\
		\hline
		$t=7$ & $q \geq 107$ and $n \geq 19$ \\
		\hline
		$t=6.3$ & $q \geq 211$ and $n \geq 13$ \\
		\hline
		$t=6.3$ & $q \geq 211$ and $n \geq 13$ \\
		\hline
		$t=6.6$ & $q \geq 980$ and $n =12$ \\
		\hline
		$t=6.8$ & $q \geq 14459$ and $n =11$ \\
		\hline
		$t=7.4$ & $q \geq 3.63\cdot 10^6$ and $n =10$ \\
		\hline
		$t=8.2$ & $q \geq 2.24\cdot 10^{13}$ and $n =9$ \\
		\hline
		$(t,u)=(6,7)$ & $q \geq 7.05\cdot 10^{21}$ and $n =8$ \\
		\hline
		$(t,u)=(8.5,9.5)$ & $q \geq 8.66 \cdot 10^{184}$ and $n =7$ \\
		\hline
	\end{tabular}\vspace*{0.5cm}
	\caption{
		Values of $q$ and  $n$ such that there exists a $3$-primitive element in  $\F_{q^n}$ which is $3$-normal over $\F_{q}$.}
	\label{table1}
\end{table}
\end{proposition}
\begin{proof}
Let $r=3$ and $k=3$. 
Using SageMath (see \cite{SAGE}) we get that condition \eqref{condition} is satisfied for 
the pairs $(q,n)$, with $n \geq 9$, given in Table \ref{table1}. For $n=8$ with $t=10$
condition \eqref{condition} is satisfied for  $q \geq 6.88\cdot 10^{51}$ and
for $n=7$ with $t=15.6$ condition \eqref{condition} is satisfied for  $q \geq 5.71\cdot 10^{3157}$.
For $n=8$ we use Lemma \ref{casetu} with
$t=6$ and $u=7$ and we get $q \geq 7.05\cdot 10^{21}$.
For $n=7$ we also use Lemma \ref{casetu} with
$t=8.5$ and $u=9.5$ and we get $q \geq 8.66 \cdot 10^{184}$.

From Proposition \ref{cota}  we get the desired result for $n \geq 9$.
and from Lemma \ref{casetu} we get the desired result for $n=7$ and $n=8$.

\end{proof}

The next lemma gives a better bound for $n=7$.

\begin{lemma}\label{casen7}
Let $q$ be a prime power such that $2.132 \cdot 10^{15} \leq q < 8.66 \cdot 10^{184}$.
If $3 \mid q^7-1$ and
$x^7-1$ has a degree $3$ factor in $\mathbb{F}_q[x]$ then
there exists a $3$-primitive element in  $\F_{q^7}$ which is $3$-normal over $\F_{q}$.
\end{lemma}
\begin{proof}
Let $q$ be a prime power such that $q < 8.66 \cdot 10^{184}$,
$3 \mid q^7-1$ and
$x^7-1$ has a degree $3$ factor $f \in \mathbb{F}_q[x]$. We will use Proposition \ref{sieve-prop} with
$\ell=q-1$, and $g=f$ if $7 \nmid q$ or $g=1$ if $7 \mid q$ (so $\widetilde{g}=\gcd(g,\frac{x^7-1}{f})=1$,
$s=s'\leq 4$ and $\delta'=0$).
Let $p$ be a prime number. If $p \mid q^7 -1$ but $p \nmid q-1$ then $7 \mid \varphi(p)=p-1$. This means
that the set $\{p_1,\ldots , p_v \}$ is composed by primes of the form $7j+1$.
Let $S_k$ and $P_k$ be, respectively, the sum of the inverses
and the product of the first $k$ primes of the form $7j+1$.
Since $\{p_1,\ldots , p_v \}$ is a set of prime numbers which divide $q^6+q^5+q^4+q^3+q^2+q+1$,
then $P_v \leq q^6+q^5+q^4+q^3+q^2+q+1 < 4.22 \cdot 10^{1109}$, 
therefore $v \leq 299$ and $S_v < 0.19113$. If we suppose $q>10^5$ then
$$
\delta = 1-\sum_{i=1}^{v}\frac{1}{p_i}-\sum_{i=1}^{s}\frac{1}{q^{\deg P_i}}
\geq 1-S_v - \frac{4}{q} > 0.80883
$$
and
$\Delta = 2 + \frac{v+s-1}{\delta} <1310.0623$.
So, observing that if $q \geq \left( 3 \cdot 1310.0623 \cdot A_t \right)^{\frac{2t}{t-2}}$ 
for some real number $t>2$, then 
$q^{\frac{7}{2}-3} \geq  3 \cdot A_t \cdot q^{\frac{1}{t}}  \cdot 1310.0623 >
3 \cdot W(q-1) \cdot  W(1) \cdot \Delta$, and using 
Proposition \ref{sieve-prop}, there exists a $3$-primitive element
in $\F_{q^7}$ which is $3$-normal over $\F_{q}$.
For $t=5.4$, the condition above becomes
$q \geq 2.132 \cdot 10^{15}$.
\end{proof}

From all the results above and using SageMath, we get the following result for finite fields of characteristic $11$.

\begin{corollary}
Let $n \geq 7$ be a positive integer and let $q=11^s$ be a power of $11$.
If $3 \mid q^n-1$ and
$x^n-1$ has a degree $3$ factor in $\mathbb{F}_q[x]$ then
there exists a $3$-primitive element in  $\F_{q^n}$ which is $3$-normal over $\F_{q}$.
\end{corollary}
\begin{proof}
Suppose first that $n \geq 8$. From Proposition \ref{caseall}
and Proposition \ref{caser3k3}, we only need to test condition
$q^{\frac{n}{2}-3} \geq 3 W(q^n-1)W(\frac{x^n-1}{f})$ for the finite number of pairs $(q=11^s,n)$ which are
not in Table \ref{table1} such that $3 \mid 11^{sn}-1$ and there exists a factor $f \in \mathbb{F}_q(x]$ of degree
$3$ of $x^n-1$. Using SageMath we get that the inequality 
$q^{\frac{n}{2}-3} \geq 3 W(q^n-1)W(\frac{x^n-1}{f})$ holds except for the pairs
$(11, 8)$,
$(11^2, 8)$,
$(11^3, 8)$,
$(11^4, 8)$,
$(11^6, 8)$,
$(11^2, 9)$,
$(11, 10)$,
$(11^2, 10)$,
$(11, 12)$,
$(11^2, 12)$.

Suppose now that $n=7$. From Proposition \ref{caseall},
Proposition \ref{caser3k3} and Lemma \ref{casen7}, we only need to test condition
$q^{\frac{7}{2}-3} \geq 3 W(q^7-1)W(\frac{x^7-1}{f})$ for $q=11^s<2.132 \cdot 10^{15}$ and $s\leq 14$ even.
Using SageMath we get that this condition 
holds only for $q=11^{14}$. We also get that condition
$q^{\frac{n}{2}-k} \geq r W(\ell) W(\gcd(g,\frac{x^n-1}{f}))\Delta$, from Proposition \ref{sieve-prop},
holds for the pairs $(11^6, 7)$, $(11^8,7)$, $(11^{10},7)$, $(11^{12},7)$,
$(11^3,8)$, $(11^4,8)$, $(11^6,8)$, $(11^2,9)$, $(11^2,10)$, $(11,12)$ and $(11^2,12)$. 
For the remaining pairs $(11^2,7)$, $(11^4,7)$, $(11,8)$,
$(11^2,8)$ and $(11,10)$,
we explicitly found a $3$-primitive, $3$-normal element.
\end{proof}

\end{document}